\documentclass[envcountsame]{mmnp}
\usepackage[tbtags]{amsmath}
\usepackage{amssymb}
\usepackage{graphicx}
\usepackage[latin1]{inputenc}

\idline{Vol. 7, No. 2}{35}
\doi{10.1051/mmnp/20127102}

\numberwithin{equation}{section}
\def\Re{\mathop{\rm Re}}
\def\Im{\mathop{\rm Im}}

\begin{document}

\title{Spectral properties of Schr\"{o}dinger-type operators and
large-time behavior of the solutions to the corresponding wave equation}


\author{A. G. Ramm\inst{1}   \thanks{\email {ramm@math.ksu.edu}} }

\vspace{0.5cm}

\institute{\inst{1} Department of Mathematics \\ Kansas State
University, Manhattan, KS 66506-2602, USA }


\abstract{Let $L$ be a linear, closed, densely defined in a Hilbert space
operator, not necessarily selfadjoint.

Consider the corresponding wave equations
\begin{align*}
    &(1) \quad \ddot{w}+ Lw=0, \quad w(0)=0,\quad \dot{w}(0)=f,
\quad \dot{w}=\frac{dw}{dt}, \quad f \in H. \\
    &(2) \quad \ddot{u}+Lu=f e^{-ikt}, \quad u(0)=0, \quad \dot{u}(0)=0,
\end{align*}
where $k>0$ is a constant.
Necessary and sufficient conditions are given for the operator $L$ not
to have eigenvalues in the half-plane Re$z<0$ and not to have a positive
eigenvalue at a given point $k_d^2 >0$.
These conditions are given in terms of the large-time
behavior of the solutions to problem (1) for generic $f$.

Sufficient conditions are given for the validity of a
version of the limiting amplitude principle for the operator $L$.

A relation between the limiting amplitude principle and the limiting
absorption principle is established. }





\titlerunning{Spectral properties of Schr\"{o}dinger-type operators}

\maketitle


\section{Introduction} \label{Intro} Let $L$ be a linear, densely defined,
closed operator in a Hilbert space $H$. Our results and techniques
are valid in a Banach space also, but we wish to think about $L$ as of a
Schr\"{o}dinger-type operator in a Hilbert space and, at times, think
that $L$ is
selfadjoint. For a Schr\"{o}dinger operator $L=-\nabla^2 + q(x)$ the
resolvent $(L-k^2)^{-1}, Im k >0$, is an integral operator with a kernel
$G(x,y,k)$ , its resolvent kernel. If $q$ is a real-valued function,
sufficiently rapidly decaying then $L$ is selfadjoint, $G(x,y,k)$ is
analytic with respect to $k$ in the half-plane Im$ k >0$, except,
possibly, for a finitely many simple poles $i k_j$, $k_j >0$,  the
semiaxis
$k \geq 0$ is filled with the points of absolutely continuous spectrum
of $L$, and there exists a
limit
$$\lim_{\epsilon \to 0}G(x,y,k+i\epsilon)=G(x,y,k)$$
for all $k>0$.

Sufficient conditions for $k^2=0$ not to be an eigenvalue of
$L$ are found in papers
\cite{R203}, \cite{R227}. Spectral analysis of the Schr\"{o}dinger
operators is
presented in many books (see, for example,  \cite{P} and \cite{S}). In
papers \cite{R4}, \cite{R11},
such an analysis was given in a class of
domains with infinite boundaries apparently for the first time,
see also \cite{R190}. In \cite{R50} an
eigenfunctions expansion theorem was proved for non-selfadjoint
Schr\"{o}dinger operators with exponentially decaying complex-valued
potential $q$. The operator $L$ in this paper is not necessarily assumed
to be selfadjoint.

In \cite{E} the validity of the limiting amplitude principle for some
class of selfadjoint operators $L$ has been established.

This principle says that, as $t \to \infty$, the solution to problem
\begin{equation} \label{eq1}
    \ddot{u} +Lu=f e^{-ikt}, \quad u(0)=0, \quad \dot{u}(0)=0,
\quad \dot{u}=\frac{du}{dt} ,
\end{equation}
has the following asymptotics
\begin{equation} \label{eq2}
    u=e^{-ikt}v +o(1), \quad t \to \infty,
\end{equation}
where $k$ is a real number and $v \in H$ solves the equation
\begin{equation} \label{eq3}
    Lv-k^2v=f.
\end{equation}
The $v$ is called {\it the limiting amplitude}. It turns out that a more
natural definition of the limiting amplitude is:
\begin{equation} \label{eq4}
    v=\lim_{t \to \infty} \frac{1}{t} \int_0^t u(s) e^{iks}ds,
\end{equation}
if this limit exists and solves equation \eqref{eq3}.

Why is this
definition more natural than \eqref{eq2}? There are good reasons for this.
One
of the reasons is: if \eqref{eq2} and \eqref{eq3} hold, then the limit
\eqref{eq4} exists and solves equation \eqref{eq3}. The other reason is:
the limit \eqref{eq4} may exist and solve equation  \eqref{eq3} although
the limit \eqref{eq2} does not exist.

{\bf Example.} {\it If $u=e^{ikt}v+e^{ik_1 t}v_1$, then the limit 
\eqref{eq2} does not exist, while the limit \eqref{eq4} does exist and is 
equal to $v$.}

To describe our assumptions and results, some preparation is needed.

Consider the problem
\begin{equation} \label{eq5}
    \ddot{w}+Lw=0,\quad w(0)=0, \quad \dot{w}(0)=f. \end{equation}
Assuming that $||u(t)|| \leq c e^{at}$, where $c >0$ stands throughout the
paper for various generic constants, and $a \geq 0$ is a constant, one
can define the Laplace transform of $u(t)$,
$$\mathcal{U}:=\mathcal{U}(p):=\int_0^\infty e^{-pt}u(t)dt,\quad \sigma
>a,$$
where $p=\sigma +i\tau$, Re$ p=\sigma$.

Let us
take the Laplace transform of \eqref{eq1} and of \eqref{eq5} to get
\begin{equation} \label{eq6}
    L\mathcal{U}+p^2\mathcal{U}=\frac{f}{p+ik},
\end{equation}
and
\begin{equation} \label{eq7}
    L\mathcal{W}+p^2\mathcal{W}=f,
\end{equation}
where
$$\mathcal{W}=\mathcal{W}(p)=\int_0^\infty w(t)e^{-pt}dt.$$
We also denote $\mathcal{W}(p):=\bar{w}(t)$.

The complex plane $p$ is related to the complex plane $k$ by the formula
\begin{equation} \label{eq8}
    p=-ik, \quad k =k_1+ik_2,\quad k_2 \geq 0, \quad \sigma=k_2\ge 0.
\end{equation}
We assume throughout that $f$ is generic in
the following sense:

{\it If $I$ is the identity operator and a point $p$ is a pole of the
kernel of the operator $(L+p^2I)^{-1}$, then it
is a pole of the same order of the element $(L+p^2I)^{-1}f=\mathcal{W}$.}

If $k^2$ is an eigenvalue of $L$ and $\Re k^2 <0$, then $\Im k >0$, where
$k=|k|e^{\frac{i \arg k^2}{2}}, p=-ik$, so $\sigma=\Re p >0$.
Let  $k>0$ and assume that $-k^2 <0$ is an eigenvalue of $L$. Then
$ik$ is a pole of the resolvent kernel $G(x,y,k)$, and $p=-i(ik)=k$ is a
pole of the kernel of the operator $(L+p^2I)^{-1}$.  If
$k^2>0$ is an eigenvalue of $L$, then $p=-ik$ is a pole of
the operator $(L+p^2I)^{-1}$.

The following known facts from the theory of Laplace transform will be used.

\begin{prpstn} \label{prop1}
An analytic in the half-plane $\sigma > \sigma_0\ge 0$ function $F(p)$ is
the
Laplace transform of a function $f(t)$, such that $f(t)=0$ for $t<0$ and
\begin{equation} \label{eq9}
    \int_0^\infty |f(t)|^2e^{-2\sigma_0 t}dt<\infty
\end{equation}
if and only if
\begin{equation} \label{eq10}
    \sup_{\sigma > \sigma_0}\int_{-\infty}^\infty |F(\sigma+i\tau)|^2d\tau
< \infty.
\end{equation}
\end{prpstn}

\begin{prpstn} \label{prop2}
If $F(p)=\overline{f(t)}$, then
\begin{equation} \label{eq11}
    \frac{F(p)}{p}=\overline{\int_0^t f(s)ds}. \end{equation} \end{prpstn}
Let us now formulate the main Assumptions A and B standing throughout this
paper.

{\bf Assumption A.} For a generic $f$
the $\mathcal{W}(p)=(L+p^2)^{-1}f$ is analytic in the half-plane $\sigma
>0$, except, possibly, at a finitely many simple poles at the
points $-ik_j$, $1 \leq j
\leq J$, $k_j$ are real numbers, and at the points $\kappa_m,$ $\Re
\kappa_m >0$,
\begin{equation} \label{eq12}
    \mathcal{W}(p)= \sum_{j=1}^J \frac{v_j}{p+ik_j}+\mathcal{W}_1(p)+
\sum_{m=1}^M \frac{b_m}{p-\kappa_m},
\end{equation}
where $v_j$ and $b_m$ are some elements of $H$, $\mathcal{W}_1(p)$
is an analytic function in the half-plane Re$p=\sigma>0$, continuous up to
the imaginary axis $\sigma=0$, and satisfying the following estimate
\begin{equation} \label{eq13}
    ||\mathcal{W}_1(p)|| \leq \frac{c}{1+|p|^\gamma}, \quad \gamma > \frac{1}{2}.
\end{equation}
{\bf Assumption B.}
There exists the limit
\begin{equation} \label{eq14}
    \lim_{\sigma \to 0} ||\mathcal{W}_1(\sigma-ik)-\mathcal{W}_1(-ik)|| =0
\end{equation}
for all real numbers $k$.
\begin{thrm} \label{thm1}
Let the Assumption A hold. Then a necessary and sufficient condition for
the
operator $L$  to have no eigenvalues in the half-plane $\Re k^2 <
0$ is the validity of the estimate
\begin{equation} \label{eq15}
    \left|\left|\int_0^t w(s)ds \right|\right|=O(e^{\epsilon t}),
\quad t\to\infty,
\end{equation}
for an arbitrary small $\epsilon >0$.

A necessary and sufficient condition for the operator $L$ not
to have any positive eigenvalues $k^2 >0$ is the validity of the estimate
\begin{equation} \label{eq16}
    \left|\left|\frac{1}{t}\int_0^t e^{iks}w(s)ds \right|\right|=o(1),
\quad t\to\infty, \quad \forall k \in \mathbb{R}.
\end{equation}
A point $ik_0 >0, k_0 >0$, is not a pole of the resolvent kernel
of the operator $(L-k^2-i0)^{-1}$ if and only if estimate  \eqref{eq16}
holds with $k=k_0>0$.
\end{thrm}

{\bf Remark.} If condition \eqref{eq16} holds for $k=0$, then
$||\int_0^t w(s)ds||=o(t)$, so condition \eqref{eq15} holds,
and the operator $L$ has no eigenvalues in the half-plane 
 $\Re k^2 <0$.

\begin{thrm} \label{thm2} Let the Assumptions A and B hold. Suppose
that estimates
\eqref{eq14} and \eqref{eq15}
hold.
Then the limiting amplitude principle \eqref{eq4} holds for every
$k\in \mathbb{R}$, $k \neq
k_j, 1 \leq j \leq J$. \end{thrm}

In section 2, proofs are given.

\section{Proofs} \subsection{Proof of Theorem 1.3} From the
Assumption A and Proposition 1.1, it follows that
$\mathcal{W}(p)$ is a Laplace transform of a function $w(t)$ such that
\begin{equation} \label{eq17}
    w(t)=\sum_{j=1}^J v_j e^{-ik_j t}+\sum_{m=1}^M b_m e^{\kappa_m t}
+w_1(t),
\end{equation}
where
\begin{equation} \label{eq18}
    w_1(t)=\frac{1}{2\pi i}\int_{\sigma_0-i\infty}^{\sigma_0+i\infty}
e^{pt} \mathcal{W}_1(p)dp, \end{equation} and the integral in \eqref{eq18}
converges in
$L^2$-sense due to the assumption \eqref{eq13}. It is clear from formula
\eqref{eq17} that all $b_m=0$ if and only if estimate \eqref{eq15} holds
with $0 < \epsilon <\min_{1\leq m \leq M}\Re \kappa_m$. This proves the
first
conclusion of Theorem 1.3.

Let us calculate the expression on the left side of formula \eqref{eq16}
and
show that this expression is $o(1)$ unless $k=k_j$ for some $1\leq j\leq
J$. In this calculation it is assumed that $L$ does not have any
eigenvalues in the half-plane $\Re k^2 < 0$, in other words, that all
$b_m=0$. Otherwise the expression on the left of formula \eqref{eq16}
tends
to infinity as $t \to \infty$ at an exponential rate.

If all $b_m=0$ in
\eqref{eq17}, then \begin{equation} \label{eq19}
    \sum_{j=1}^J v_j \frac{1}{t} \int_0^t e^{i(k-k_j)t}dt +\frac{1}{t}
\int_0^t w_1(t) e^{ikt}dt := I_1+I_2.
\end{equation}
If $k$ and $k_j$ are real numbers, then
\begin{equation} \label{eq20}
    \lim_{t\to \infty} \frac{1}{t} \int_0^t e^{i(k-k_j)t} dt =
    \left\{
        \begin{array}{ll}
            1, & \quad k=k_j, \\
            0, & \quad k \neq k_j.
        \end{array}
    \right.
\end{equation}
Thus, $I_1=0$ if and only if $k$ does not coincide with any of $k_j$,
$1\le j \le J$.

Let us prove that
\begin{equation} \label{eq21}
    \lim_{t\to \infty} \frac{1}{t} \int_0^t w_1(t)e^{ikt} dt = 0.
\end{equation}
By proposition \eqref{prop2} and the Mellin inversion formula, one has
\begin{equation} \label{eq22}
    I:= \frac{1}{t} \int_0^t w_1(t)e^{ikt} dt = \frac{1}{2 \pi i}
\int_{\sigma-i\infty}^{\sigma+i\infty} \mathcal{W}_1(p-ik)\frac{e^{pt}}{pt}dp,
\end{equation}
where Re$p=\sigma >0$ can be chosen arbitrarily small.

Let $pt=q$,  take $\sigma=\frac{1}{t}$, write $q=1+is$, and write the
integral on the right side of \eqref{eq22} as:
\begin{equation} \label{eq23}
    I= \frac{1}{2 \pi i}\int_{1-i\infty}^{1+i\infty}
\mathcal{W}_1(\frac{q}{t}-ik)\frac{q}{t} \frac{e^{q}}{q^2}dq.
\end{equation}
If one uses estimate \eqref{eq13} and formula $|q|=(1+s^2)^{1/2}$, then
one obtains the following inequality
\begin{equation} \label{eq24}
    ||I|| \leq \frac{1}{2\pi t}\int_{-\infty}^{\infty} \frac{1}{(1+s^2)^{1/2}}
\frac{c ds}{[1+|\frac{1+is}{t}-ik|^\gamma]}= \frac{c}{2 \pi t^{1-\gamma}}
\int_{-\infty}^{\infty}\frac{1}{(1+s^2)^{1/2}}\frac{ds}{(t^\gamma +
[1+(s-kt)^2]^{\gamma/2})}.
\end{equation}
Let $s=ty$. Then the integral on the right side of \eqref{eq24} can be written
as
\begin{align} \label{eq25}
    &\frac{ct}{2 \pi t^{1-\gamma}}\int_{-\infty}^{\infty}\frac{dy}
{(1+t^2y^2)^{1/2}}\frac{1}{(t^\gamma +[1+t^2(y-k)^2]^{\gamma/2})}
\nonumber \\
    &=\frac{c}{2 \pi}\int_{-\infty}^{\infty}\frac{dy}{(1+t^2y^2)^{1/2}}
\frac{1}{(1 +[t^{-2}+(y-k)^2]^{\gamma/2})} \nonumber \\
    &\leq \frac{c}{2 \pi}\int_{-\infty}^{\infty}\frac{dy}{(1+t^2y^2)^{1/2}}
\frac{1}{[1 +(y-k)^{\gamma}]} \to 0,
\text{ as } t \to \infty,
\end{align}
and the convergence of the last integral to zero is uniform with respect
to $k\in \mathbb{R}$.

Thus
\begin{equation} \label{eq26}
    \lim_{t \to \infty}||I||=0.
\end{equation}
From \eqref{eq19}-\eqref{eq21} the last two conclusions of
Theorem 1.3 follow. Theorem 1.3 is proved.
\hfill$\Box$

\subsection{Proof of Theorem 1.4} Using Proposition 1.2,
 equation \eqref{eq6}, and the Mellin formula, one gets
\begin{equation}
\label{eq27}
    \frac{1}{t} \int_0^t u(t)e^{ikt} dt = \frac{1}{t}\frac{1}{2 \pi i}
\int_{\sigma-i\infty}^{\sigma+i\infty} \frac{\mathcal{U}(p-ik)}{p}e^{pt}dp,
\end{equation}
where, according to \eqref{eq6},
\begin{equation} \label{eq28}
    \mathcal{U}(p-ik)=\frac{\mathcal{W}(p-ik)}{p}.
\end{equation}
Let $\sigma=\frac{1}{t}$ and $pt=q$. Then
\begin{equation} \label{eq29}
    \frac{1}{t} \int_0^t u(t)e^{ikt} dt = \frac{1}{2 \pi i}
\int_{1-i\infty}^{1+i\infty} \mathcal{W}\left(\frac{q}{t}-ik\right)
\frac{e^q}{q^2}dq. \end{equation} Estimate \eqref{eq15} and Theorem 1.3
imply that all
$b_m=0$ in formula \eqref{eq17}.
Therefore,
using formula \eqref{eq17} with $b_m=0$, one gets
$$\mathcal{W}=\sum_{j=1}^{J}v_j \frac{1}{p+ik_j} +\mathcal{W}_1,$$
and
\begin{equation}
\label{eq30}
    \mathcal{W}\left(\frac{q}{t}-ik\right)=\mathcal{W}_1\left(\frac{q}{t}-
ik\right)+\sum_{j=1}^{J}v_j \frac{1}{\frac{q}{t}-i(k-k_j)}.
\end{equation}
One has $\overline{t^n}=\frac{n!}{p^{n+1}}$. Therefore
\begin{equation*}
     \frac{1}{2 \pi i}\int_{1-i\infty}^{1+i\infty} \frac{e^q}{q^2}dq=1,
\end{equation*}
and
\begin{equation} \label{eq31}
    \lim_{t \to \infty} \frac{v_j}{2 \pi i}\int_{1-i\infty}^{1+i\infty}
\frac{1}{\frac{q}{t}-i(k-k_j)}\frac{e^q}{q^2}dq=
    \left\{
        \begin{array}{ll}
            \frac{iv_j}{k-k_j},& \quad k\neq k_j, \\
            \infty ,& \quad k= k_j.
        \end{array}
    \right.
\end{equation}
Furthermore,
\begin{equation} \label{eq32}
    \lim_{t \to \infty} \frac{1}{2 \pi i}\int_{1-i\infty}^{1+i\infty}
\mathcal{W}_1\left(\frac{q}{t}-ik\right)\frac{e^q}{q^2}dq=\mathcal{W}_1(-ik),
\end{equation} as follows from assumption \eqref{eq14} and the
Lebesgue's dominated convergence theorem if one passes to the limit $t
\to \infty$
under the sign of the integral \eqref{eq32}.
 Let us check that this $v$ solves equation \eqref{eq3}. This would
conclude the proof of Theorem 1.4. We need a lemma.

{\bf Lemma 1.} {\it If $h\in L^1_{loc}(0,\infty)$and the limit
$\lim_{t\to \infty}t^{-1}\int_0^th(s)ds$ exists, then the limit
$\lim_{p\to 0}p\int_0^{\infty}e^{-pt}h(t)dt$ exists, and
\begin{equation} \label{eq33}
\lim_{t\to \infty}t^{-1}\int_0^t h(s)ds=\lim_{p\to
0}p\int_0^{\infty}e^{-pt}h(t)dt.
\end{equation}
} {\bf Proof of Lemma 1}. One has
$$p\int_0^{\infty}e^{-pt}h(t)dt=pe^{-pt}\int_0^t
h(s)ds|_{0}^{\infty} +p^2\int_0^{\infty}te^{-pt}t^{-1}\int_0^t h(s)ds dt.
$$
For any $p>0$ one has
$$pe^{-pt}\int_0^t h(s)ds|_{0}^{\infty}=0.$$
Let $q=pt$ and denote $H(t):=t^{-1}\int_0^t h(s)ds$, $J:=\lim_{t\to
\infty}H(t)$. Then
$$\lim_{p\to 0}p^2\int_0^{\infty}te^{-pt}t^{-1}\int_0^t h(s)ds dt=
\lim_{p\to 0}\int_0^{\infty}qe^{-q}H(qp^{-1})dq.$$
Passing in the last integral to the limit $p\to 0$ one obtains
\eqref{eq33}. Lemma 1 is proved. \hfill $\Box$

Using equation \eqref{eq33}, one writes $v=\lim_{p\to
0}p\mathcal{U}(p-ik)$, where $\mathcal{U}$ solves equation
\eqref{eq6}. Thus,
$$L\mathcal{U}(p-ik)+(p-ik)^2\mathcal{U}(p-ik)=p^{-1}f.$$
Multiplying both sides of this equation by $p$ and passing to the
limit $p\to 0$, one obtains equation \eqref{eq3}. In the passage
to the limit under the sign of the unbounded operator $L$
the assumption that $L$ is closed was used.

Thus,  the conclusion of Theorem 1.4 follows. \hfill$\Box$

If the limit \eqref{eq14} exists at a point $p=i\tau$ then one says that
the limiting absorption principle holds for the operator $L$ at the point
$k=ip=i(-ik)=k, k>0$.

Thus, Assumption B means that the limiting absorption principle holds for
$L$ at the point $k>0$, that is, $\lim_{\epsilon \to
0}(L-k^2-i\epsilon)^{-1}f$ exists.

\section{Applications} Let $L=-\nabla^2+q(x)$, where $q(x)$ is a
real-valued function, $|q(x)|\leq c(1+|x|)^{-2-\epsilon}, \epsilon >0, x
\in \mathbb{R}^3$. Then $L$ is selfadjoint on the domain
$H^2(\mathbb{R}^3)$. Its
resolvent $(L-k^2-i0)^{-1}$ satisfies Assumptions A and B if one keeps in
mind the following.

Let $G(x,y,k)$ be the resolvent kernel of $L$, that is, the kernel of the
operator $(L-k^2-i0)^{-1}$,
$$LG(x,y,k)=-\delta(x-y) \quad  in \quad \mathbb{R}^3,$$
$G\in L^2(\mathbb{R}^3)$ for $\Im k>0$. If $f \in L^2(\mathbb{R}^3)$ is
compactly supported, then for $k>0$ the function
$$v(x):=(L-k^2-i0)^{-1}f=\int_{\mathbb{R}^3} G(x,y,k) f(y)dy$$ does not
necessarily belong
to $L^2(\mathbb{R}^3)$.

For example, if $q(x)=0$, then $G(x,y,k)=\frac{e^{ik|x-y|}}{4\pi |x-y|}$,
and the function \begin{equation} \label{eq34}
     v(x,k)=\int_{|y|\leq 1} g(x,y,k)dy =O\left(\frac{1}{|x|}\right)
\end{equation}
does not belong to $L^2(\mathbb{R}^3)$ (except for those $k>$ for which
$x(x,k)=0$ in the region $|y|\ge 1$. These numbers $k>0$ are
the zeros of the Fourier transform of the characteristic function of
the ball $|y|\le 1$,  see \cite{R470}, Chapter 11.

By this reason the abstract results of theorem \eqref{thm1} and
\eqref{thm2} can be used in applications if one defines some subspace of
$H$, for example, a subspace of functions with compact support, denote by
$\mathcal{P}$, a projection operator on this subspace, and replaces
$\mathcal{W}$ and $\mathcal{W}_1$ by $\mathcal{PW}$ and $\mathcal{PW}_1$
in equations \eqref{eq12} and  \eqref{eq14}. For example, the function 
\eqref{eq34} one
replaces by $\eta(x)v(x,k)$, where $\eta(x)$ is a characteristic function
of a compact subset of $\mathbb{R}^3$.

The analytic properties of $\eta(x)v(x,k)$ and of $v(x,k)$ as functions of
$k$ are the same. A similar suggestion is used in \cite{E}.

With the above in mind, one knows (for example, from \cite{P} or \cite{S})
that Assumptions A and B hold for $L=-\nabla^2+q(x)$.

Consequently, the conclusions of Theorems 1.3 and 1.4
hold.

In addition, the assumptions
$$|q(x)|\leq c(1+|x|)^{-2-\epsilon},\quad
\epsilon >0, \quad \Im q =0,$$ imply that $L$ does not have positive
eigenvalues,
so all $v_j=0$, and zero is not an eigenvalue of $L\geq 0$ if $\epsilon
>0$ (see \cite{R203}, \cite{R227}).

A new method for estimating of large time behavior of solutions to
abstract evolution problems is developed in \cite{R605}, where some
applications of this method are given.





\end{document}